\documentclass{amsart}

\usepackage{amssymb}
\usepackage{amsmath}
\usepackage[all]{xy}
\usepackage{enumerate}
\usepackage[usenames,dvipsnames]{pstricks}
\usepackage{epsfig}
\usepackage{pst-grad} 
\usepackage{pst-plot} 

\newtheorem{theorem}[subsection]{Theorem}

\theoremstyle{remark}
\newtheorem{exa}[subsection]{Example}
\newtheorem{rem}[subsection]{Remark}

\theoremstyle{definition}

\DeclareMathOperator{\Int}{Int}

\DeclareMathOperator{\codim}{codim}
\DeclareMathOperator{\BOX}{Box}
\DeclareMathOperator{\Hom}{Hom}

\DeclareMathOperator{\Span}{span}

\DeclareMathOperator{\orb}{orb}

\DeclareMathOperator{\bd}{bd}
\DeclareMathOperator{\age}{age}

\begin{document}

\title{Ehrhart Theory for Lawrence Polytopes and Orbifold Cohomology of Hypertoric Varieties}

\author{Alan Stapledon}
\address{Department of Mathematics, University of Michigan, Ann Arbor, MI 48109, U.S.A.}
\email{astapldn@umich.edu}


\thanks{The author would like to thank Nicholas Proudfoot 
for some useful comments.
}


\begin{abstract}
We establish a connection between the orbifold cohomology of hypertoric varieties and the Ehrhart theory of Lawrence polytopes.
More specifically, we show that the dimensions of the orbifold cohomology groups of a hypertoric variety are equal to the coefficients of the Ehrhart $\delta$-polynomial of the associated Lawrence polytope. As a consequence, we deduce a formula for the Ehrhart $\delta$-polynomial of a Lawrence polytope and use the injective part of the Hard Lefschetz Theorem for hypertoric varieties to deduce some inequalities between the coefficients of the $\delta$-polynomial.
\end{abstract}

\maketitle




\section{Introduction}

Hypertoric varieties were introduced by Biewlawski and Dancer \cite{BDGeometry} using a hyperk\"ahler analogue of the construction of toric varieties
by K\"ahler quotients. Hausel and Sturmfels \cite{HauQuaternionic,HSToric} showed that the cohomology of hypertoric varieties is intimately related to the combinatorics of matroids and hyperplane arrangements. Recently, Jiang and Tseng \cite{JTOrbifold} and Goldin and Harada \cite{GHOrbifold} independently gave explicit descriptions of the orbifold cohomology rings of hypertoric varieties. We will establish a relationship between the orbifold cohomology of  hypertoric varieties and the lattice point enumeration of Lawrence polytopes.

Let $N$ be a lattice of rank $d$ and let $\mathcal{B} = \{ b_{1}, \ldots, b_{n} \}$ be a configuration of $n$ non-zero lattice points in $N$  that  generate $N$ as an abelian group.
Let $M = \Hom_{\mathbb{Z}}(N,\mathbb{Z})$ denote the dual lattice to $N$ and let $\mathcal{R} = \{r_{1}, \ldots, r_{n} \}$ be a collection of real numbers.
Consider the \emph{weighted, co-oriented hyperplane arrangement} $\mathcal{H}$ in $M_{\mathbb{R}} = M \otimes_{\mathbb{Z}} \mathbb{R}$ defined by 
the hyperplanes
$H_{i} = \{ u \in M_{\mathbb{R}} \mid \langle u, b_{i} \rangle = r_{i}     \}$,
for $i = 1, \ldots, n$, and the co-normal vectors $\mathcal{B}$. We may and will choose $\mathcal{R}$ so that the hyperplane arrangement $\mathcal{H}$ is \emph{simple}. That is, the intersection of any $l$ hyperplanes in $\mathcal{H}$ is either empty or has codimension $l$.
The \emph{Lawrence toric variety} $X_{\mathcal{H}}$ 
associated to $\mathcal{H}$ is a $(d + n)$-dimensional simplicial toric variety,
 and the corresponding hypertoric variety $\mathcal{M}_{\mathcal{H}}$ 
 is a $2d$-dimensional complete intersection in $X_{\mathcal{H}}$, with a canonical orbifold structure \cite[Section 6]{HSToric}. In Section \ref{construction}, we will give a concrete local construction of $X_{\mathcal{H}}$ and $\mathcal{M}_{\mathcal{H}}$ in terms of $\mathcal{H}$, expanding on the descriptions in \cite{HSToric} and \cite{JTOrbifold}. We refer the reader to \cite{HSToric} and \cite{ProSurvey} for an algebraic description of hypertoric varieties using Gale diagrams and for the construction of hypertoric varieties as hyperk\"ahler quotients.
 
 The theory of orbifold cohomology, developed by Chen and Ruan \cite{CRNew, CROrbifold},
associates to an orbifold $Y$ a finite-dimensional $\mathbb{Q}$-algebra
$H_{\orb}^{*}(Y, \mathbb{Q})$, graded by $\mathbb{Q}$. In our case, the orbifold cohomology groups of $X_{\mathcal{H}}$ and $\mathcal{M}_{\mathcal{H}}$ are graded over $\mathbb{Z}$ and depend only on $\mathcal{B}$ and not on $\mathcal{H}$ \cite[Theorem 1.1]{JTOrbifold}.

\begin{theorem}\cite[Theorem 3.10]{JTNote}\label{lemma1}
The orbifold cohomology groups $H_{\orb}^{*}(X_{\mathcal{H}}, \mathbb{Q})$ and
$H_{\orb}^{*}(\mathcal{M}_{\mathcal{H}}, \mathbb{Q})$ 
are isomorphic as graded $\mathbb{Q}$-algebras and vanish in all degrees not in $2\mathbb{Z}$.
\end{theorem}


If $e_{1}, \ldots, e_{n}$ denote the standard generators of
$\mathbb{Z}^{n}$, then the \emph{Lawrence polytope} $P_{\mathcal{B}}$ associated to $\mathcal{B}$ is the
$(d + n - 1)$-dimensional lattice polytope in $N \times \mathbb{Z}^{n}$ with 
vertices
$\{ (b_{i},e_{i}),(0,e_{i})\mid  i = 1, \ldots, n \}$.
Lawrence polytopes have been crucial in the construction of non-rational polytopes in the work of  Perles \cite{GruConvex}, Mn\"ev \cite{MnevUniversality} and
Ziegler \cite{ZieNon, ZieLectures}, and their combinatorial properties have been studied by Bayer and Sturmfels \cite{BSLawrence} and Santos \cite{SanTriangulations}. Recently, they appeared in Batyrev and Nill's classification of degree $1$ lattice polytopes \cite{BNMultiples} and in
the construction  of hypertoric varieties \cite{HSToric}.

For every positive integer $m$, let $f_{\mathcal{B}}(m) := \# \left( mP_{\mathcal{B}} \cap (N \times \mathbb{Z}^{n}) \right)$ denote the
number of lattice points in the $m$'th dilate of $P_{\mathcal{B}}$. A famous theorem of Ehrhart \cite{ehrhartpolynomial} asserts that
$f_{\mathcal{B}}(m)$ is a polynomial in $m$ of degree $\dim P_{\mathcal{B}} = n + d - 1$, called the
\emph{Ehrhart polynomial} of $P_{\mathcal{B}}$.
Equivalently,
the generating series of $f_{\mathcal{B}}(m)$ can be written in the form
\begin{equation*}
\sum_{m \geq 0} f_{\mathcal{B}}(m) \, t^{m} = \frac{ \delta_{\mathcal{B}}(t) }{ (1 - t)^{n + d} } \, ,
\end{equation*}
where $\delta_{\mathcal{B}}(t) = \delta_{0} + \delta_{1}t + \cdots +  \delta_{n + d - 1}t^{n + d - 1}$ is a polynomial of degree at most $d$ with
integer coefficients,
called the \emph{Ehrhart $\delta$-polynomial} of $P_{\mathcal{B}}$.
Ehrhart $\delta$-polynomials of lattice polytopes have been studied extensively over the last forty years by many authors including Stanley \cite{StaHilbert2, StaDecompositions, StaHilbert1,  StaMonotonicity} and Hibi \cite{HibSome, HibEhrhart,  HibLower, HibStar}.
In a recent paper \cite{YoWeightI}, the author expressed the coefficients of  the Ehrhart $\delta$-polynomial of a lattice polytope as sums of dimensions of orbifold cohomology groups of a toric orbifold. We will use the following result.

\begin{theorem}\cite[Theorem 4.6]{YoWeightI}\label{oldie}
Let $P$ be an $m$-dimensional lattice polytope in a lattice $N'$ and let $\mathcal{T}$ be a lattice triangulation of $P$.
If $\sigma$ denotes the cone over $P \times \{ 1 \}$ in $(N' \times \mathbb{Z})_{\mathbb{R}}$ and if $\triangle$ denotes the simplicial fan refinement of $\sigma$ induced by $\mathcal{T}$, then we may consider the corresponding toric variety $Y = Y(\triangle)$, with its canonical orbifold structure.
The Ehrhart $\delta$-polynomial of $P$ has the form
\[
\delta_{P}(t) = \sum_{i = 0}^{m} \dim_{\mathbb{Q}} H_{\orb}^{2i}(Y, \mathbb{Q}) t^{i}.
\]
\end{theorem}

We deduce the following geometric interpretation of $\delta_{\mathcal{B}}(t)$.

\begin{theorem}\label{punter}
The Ehrhart $\delta$-polynomial of the Lawrence polytope $P_{\mathcal{B}}$ has the form
\[
\delta_{\mathcal{B}}(t) = \sum_{i = 0}^{n + d - 1} \dim_{\mathbb{Q}} H_{\orb}^{2i}(X_{\mathcal{H}}, \mathbb{Q}) t^{i}
= \sum_{i = 0}^{n + d - 1} \dim_{\mathbb{Q}} H_{\orb}^{2i}(\mathcal{M}_{\mathcal{H}}, \mathbb{Q}) t^{i}.
\]
\end{theorem}
\begin{proof}
If $\psi: N \times \mathbb{Z}^{n} \rightarrow \mathbb{Z}$ is given by $\psi(\lambda + \sum_{i = 1}^{n} \mu_{i}e_{i} ) = \sum_{i = 1}^{n} \mu_{i}$, for $\lambda \in N$ and integers $\mu_{i}$,  then we may choose an isomorphism $N \cong \psi^{-1}(0) \times \mathbb{Z}$ such that $P_{\mathcal{B}}$
is a $(n + d - 1)$-dimensional lattice polytope in $\psi^{-1}(0) \times \{1\}$.
If $\Sigma_{\mathcal{H}}$ denotes the fan associated to the toric variety $X_{\mathcal{H}}$, then
the hyperplane arrangement $\mathcal{H}$ induces a regular lattice triangulation $\mathcal{T}$ of $P_{\mathcal{B}}$ such that the cones of $\Sigma_{\mathcal{H}}$ are given by the cones over the faces of $\mathcal{T}$ \cite[Proposition 4.2]{HSToric}.
The result now follows from
Theorem \ref{lemma1} and Theorem \ref{oldie}.
\end{proof}

Theorem \ref{punter} and computations of the orbifold cohomology ring of $\mathcal{M}_{\mathcal{H}}$ by Jiang and Tseng \cite{JTNote,JTOrbifold}  give a combinatorial formula for $\delta_{\mathcal{B}}(t)$. In order to state the result,
first recall that  the
\emph{matroid} $M_{\mathcal{B}}$ associated to  $\mathcal{B}$ is the collection of all linearly independent subsets of $\mathcal{B}$, as well as the origin $\{ 0 \}$.
The dimension of an element $F$ in $M_{\mathcal{B}}$ is equal to the dimension of $\Span F$, the 
span
of the elements of $F$ in $N_{\mathbb{R}}$. 
For any $F$ in $M_{\mathcal{B}}$, consider the projection
$\phi_{F}: N \rightarrow N/ (N \cap \Span F)$ and the induced configuration
$\mathcal{B}_{F} = \{ \phi_{F}(b_{i}) \mid b_{i} \notin \Span F\}$ in $N/ (N \cap \Span F)$, with associated matroid $M_{F}$. The hyperplane arrangement $\mathcal{H}$ restricts to a simple hyperplane arrangement $\mathcal{H}^{F}$ on the affine space $\cap_{b_{i} \in F} H_{i}$,
which we identify with $M_{F,\mathbb{R}} = \Hom_{\mathbb{R}}(N_{\mathbb{R}}/\Span F,\mathbb{R})$ after translation.
The $h$-vector of $M_{F}$ is defined by
\[
h_{F}(t) = \sum_{i = 0}^{\codim F} f_{i} t^{i} (1 - t)^{\codim F - i},
\]
where $f_{i}$ equals the number of elements in $M_{F}$ of dimension  $i$.
The hyperplane arrangement $\mathcal{H}^{F}$ divides $M_{F,\mathbb{R}}$ into locally closed cells (see Section \ref{construction})
and we will consider the $h$-vector associated to the complex of bounded cells,
\[
h^{\bd}_{F}(t) = \sum_{i = 0}^{\codim F} f^{\bd}_{i}(t - 1)^{i},
\]
where $f^{\bd}_{i}$ equals the number of bounded cells in $\mathcal{H}^{F}$ of dimension $i$.
Since $\mathcal{H}^{F}$ is simple, it follows from results of Zaslavsky that $h_{F}(t) = h^{\bd}_{F}(t)$ \cite{ZasFacing}.
The following geometric interpretation of $h_{F}(t) = h^{\bd}_{F}(t)$
was observed by Hausel and Sturmfels, after Konno had computed the cohomology ring of a hypertoric variety in \cite{KonCohomology}.

\begin{theorem}\cite[Theorem 1.1,Corollary 6.6]{HSToric}\label{hussey}
With the notation above, for any element $F$ in $M_{\mathcal{B}}$, the varieties $X_{\mathcal{H}^{F}}$ and $\mathcal{M}_{\mathcal{H}^{F}}$ have no odd cohomology and
\[ h_{F}(t) = h^{\bd}_{F}(t) = \sum_{i = 0}^{\codim F} \dim_{\mathbb{Q}} H^{2i}(X_{\mathcal{H}^{F}}, \mathbb{Q}) t^{i}
= \sum_{i = 0}^{\codim F} \dim_{\mathbb{Q}} H^{2i}(\mathcal{M}_{\mathcal{H}^{F}}, \mathbb{Q}) t^{i}. \]
\end{theorem}


If $F = \{ b_{k_{1}},  \ldots, b_{k_{r}} \}$ lies in $M_{\mathcal{B}}$, then let $\BOX(F) = \{ v \in N_{\mathbb{R}} \mid v = \sum_{j = 1}^{r} \alpha_{j} b_{k_{j}}, 0 < \alpha_{j} < 1 \}$ and let $\BOX( \{ 0 \} )  = \{ 0 \}$. The following result follows from Theorem \ref{hussey} and the work of Jiang and Tseng.

\begin{theorem}\cite[Proposition 4.7]{JTOrbifold}\label{bussey}
With the notation above,
\[
\sum_{i = 0}^{n + d - 1} \dim_{\mathbb{Q}} H_{\orb}^{2i}(X_{\mathcal{H}}, \mathbb{Q}) t^{i}
= \sum_{i = 0}^{d + n - 1} \dim_{\mathbb{Q}} H_{\orb}^{2i}(\mathcal{M}_{\mathcal{H}}, \mathbb{Q}) t^{i} = \sum_{F \in M} \# \left( \BOX(F) \cap N \right) t^{\dim F} h_{F}(t). 
\]
\end{theorem}

We deduce the following result and refer the reader
to Section \ref{combinatorics} for a combinatorial proof.

\begin{theorem}\label{newbe}
The Ehrhart $\delta$-polynomial of the Lawrence polytope $P_{\mathcal{B}}$ is equal to
\[
\delta_{\mathcal{B}}(t) = \sum_{F \in M} \# \left( \BOX(F) \cap N \right) t^{\dim F} h_{F}(t) =
\sum_{F \in M} \# \left( \BOX(F) \cap N \right) t^{\dim F} h^{bd}_{F}(t).
\]
In particular, if we write $b_{i} = a_{i}v_{i}$, for some primitive integer vector $v_{i}$ and some positive integer $a_{i}$,
and set $R^{bd}_{F}$ to be the number of bounded regions of $M_{F, \mathbb{R}} \smallsetminus \mathcal{H}^{F}$,
then
$\delta_{0} = 1$, $\delta_{1} = \sum_{i = 1}^{n}a_{i} - d$, $\delta_{d} = \sum_{F \in M} \# \left( \BOX(F) \cap N \right) R^{bd}_{F}$ and
$\delta_{d + 1} = \ldots = \delta_{d + n - 1} = 0$. If $V_{F}$ denotes the number of vertices in $\mathcal{H}^{F}$, which equals the number of maximal elements of $M_{F}$, then $(d + n - 1)!$ times the volume of $P_{\mathcal{B}}$ is equal to $\sum_{F \in M} \# \left( \BOX(F) \cap N \right) V_{F}$.
\end{theorem}

\begin{rem}
If we change the co-orientation of $\mathcal{B}$, i.e. replace some of the $b_{i}$ by $-b_{i}$, then we obtain an isomorphic hypertoric variety $\mathcal{M}_{\mathcal{H}}$ \cite[Theorem 2.2]{HPProperties}, and hence Theorem \ref{punter} implies that $\delta_{\mathcal{B}}(t)$ remains unchanged.
This can also be deduced using the formula for $\delta_{\mathcal{B}}(t)$ in Theorem \ref{newbe}.
\end{rem}

The matroid $M_{\mathcal{B}}$ is called \emph{coloop free} if there does not exist an element $b_{i}$ which lies in every maximal linearly independent subset of $\mathcal{B}$. One verifies that this condition holds if and only if $M_{\mathbb{R}} \smallsetminus \mathcal{H}$ contains a bounded region.
Hausel and Sturmfels showed that if $M_{\mathcal{B}}$ is coloop free then the
injective part of the Hard Lefschetz theorem holds for the hypertoric variety $\mathcal{M}_{\mathcal{H}}$ and hence the $h$-vector of $M_{\mathcal{B}}$ is a
\emph{MacCauley vector} \cite[Theorem 7.4]{HSToric} (c.f. \cite{SwaElements}). That is, if we write $h_{\{0\}}(t) = \sum_{i = 0}^{d} h_{i} t^{i}$, then there exists a graded $\mathbb{Q}$-algebra $R = R_{0} \oplus R_{1} \oplus \cdots
\oplus R_{\lfloor \frac{d}{2} \rfloor}$ generated by $R_{1}$ and with $g_{i} = h_{i} - h_{i -1} = \dim_{\mathbb{Q}} R_{i}$. We will use the following corollary and refer the reader to \cite{HauQuaternionic} for the definition of the $g$-inequalities.

\begin{theorem}\cite[Corollary 4.2]{HauQuaternionic}\label{rubbish}
If $M_{\mathcal{B}}$ is coloop free then its $h$-vector satisfies
$h_{i} \leq h_{j}$,
for $i \leq j \leq d - i$, and satisfies the $g$-inequalities.
\end{theorem}

We deduce the following inequalities between the coefficients of the Ehrhart $\delta$-polynomial of $P_{\mathcal{B}}$.

\begin{theorem}
If $M_{\mathcal{B}}$ is coloop free, then the Ehrhart $\delta$-polynomial of $P_{\mathcal{B}}$ satisfies
$\delta_{i} \leq \delta_{j}$, for $i \leq j \leq d - i$.
\end{theorem}
\begin{proof}
One verifies that if $M_{\mathcal{B}}$ is coloop free then $M_{F}$ is coloop free for any element $F$ in $M_{\mathcal{B}}$.
By Theorem \ref{newbe}, $\delta_{\mathcal{B}}(t) = \sum_{F \in M} \# \left( \BOX(F) \cap N \right) t^{\dim F} h_{F}(t)$.
If we write $h_{F}(t) = \sum_{i = 0}^{\codim F} h_{F,i} t^{i}$ then Theorem \ref{rubbish} implies that
$h_{F, i - \dim F} \leq h_{F, j - \dim F}$ for $i \leq j \leq d + \dim F - i$ and the result follows.
\end{proof}

\begin{rem}
The theorem above fails if $M_{\mathcal{B}}$ is not coloop free. For example, if $N = \mathbb{Z}$ and $\mathcal{B} = \{ 1 \}$, then $P$ is an interval of length $1$
and $\delta_{0} = 1 > \delta_{1} = 0$.
\end{rem}

\begin{exa}\label{visit}
Let $N = \mathbb{Z}^{2}$ and set $\mathcal{B} = \{ b_{1}, b_{2}, b_{3}, b_{4} \} = \{ (1,0), (0,1), (-2,0), (2,-1)  \}$ and
$\{ r_{1}, r_{2}, r_{3}, r_{4} \} = \{ 0, 0, 2, -1 \}$. The matroid $M_{\mathcal{B}}$ is given by $\{ 0, b_{1}, b_{2}, b_{3}, b_{4}, b_{1}b_{2}, b_{1}b_{4},
b_{2}b_{3}, b_{2}b_{4}, b_{3}b_{4} \}$ and $h_{\{ 0\}}(t) = 1 + 2t + 2t^{2}$. We have $\BOX(\{b_{2}b_{4}\}) \cap N = \{(1,0)\}$, $M_{\{b_{2}b_{4}\}} =
\{ 0 \}$ and $h_{\{b_{2}b_{4}\}}(t) = 1$. Also $\BOX(\{b_{3}\}) \cap N = \{(-1,0)\}$, $M_{\{b_{3}\}} =
\{ 0, b_{2}, b_{4} \}$ and $h_{\{b_{3}\}}(t) = 1 + t$. We conclude that $P_{\mathcal{B}}$ is a $5$-dimensional lattice polytope with $8$ vertices and
$\delta_{P}(t) = (1 + 2t + 2t^{2}) + t^{2} + t(1 + t) = 1 + 3t + 4t^{2}$.

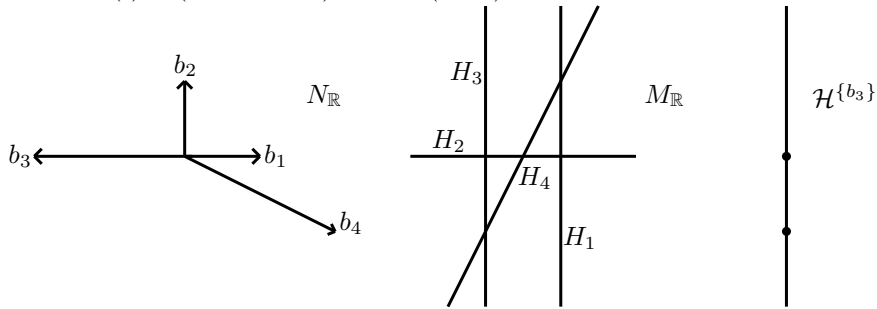
\scalebox{1} 
{
\begin{pspicture}(0,-2.02)(13.702812,2.02)
\psline[linewidth=0.04cm](7.7809377,2.0)(7.7809377,-2.0)
\psline[linewidth=0.04cm](6.7809377,2.0)(6.7809377,-2.0)
\psline[linewidth=0.04cm](6.2809377,-2.0)(8.280937,2.0)
\psline[linewidth=0.04cm](5.7809377,0.0)(8.780937,0.0)
\usefont{T1}{ptm}{m}{n}
\rput(9.162344,0.81){$M_{\mathbb{R}}$}
\usefont{T1}{ptm}{m}{n}
\rput(8.032344,-1.09){$H_{1}$}
\usefont{T1}{ptm}{m}{n}
\rput(6.532344,1.11){$H_{3}$}
\usefont{T1}{ptm}{m}{n}
\rput(6.3323436,0.21){$H_{2}$ }
\usefont{T1}{ptm}{m}{n}
\rput(7.432344,-0.29){$H_{4}$}
\psline[linewidth=0.04cm](2.7809374,0.0)(4.7809377,-1.0)
\psline[linewidth=0.04cm](4.7309375,-0.9)(4.7809377,-1.0)
\psline[linewidth=0.04cm](4.6809373,-1.04)(4.7809377,-1.0)
\psline[linewidth=0.04cm](2.7809374,0.0)(3.7809374,0.0)
\psline[linewidth=0.04cm](2.7809374,1.0)(2.7809374,0.0)
\psline[linewidth=0.04cm](0.7809375,0.0)(2.8809376,0.0)
\psline[linewidth=0.04cm](3.7809374,0.0)(3.6809375,-0.1)
\usefont{T1}{ptm}{m}{n}
\rput(4.632344,0.81){$N_{\mathbb{R}}$}
\usefont{T1}{ptm}{m}{n}
\rput(3.9923437,0.01){$b_{1}$}
\usefont{T1}{ptm}{m}{n}
\rput(2.7923439,1.21){$b_{2}$}
\usefont{T1}{ptm}{m}{n}
\rput(0.59234375,0.01){$b_{3}$}
\usefont{T1}{ptm}{m}{n}
\rput(4.992344,-0.89){$b_{4}$}
\psline[linewidth=0.04cm](3.6809375,0.1)(3.7809374,0.0)
\psline[linewidth=0.04cm](2.7809374,1.0)(2.6809375,0.9)
\psline[linewidth=0.04cm](2.7809374,1.0)(2.8809376,0.9)
\psline[linewidth=0.04cm](0.7809375,0.0)(0.8809375,0.1)
\psline[linewidth=0.04cm](0.7809375,0.0)(0.8809375,-0.1)
\psline[linewidth=0.04cm](10.780937,2.0)(10.780937,-2.0)
\psdots[dotsize=0.12](10.780937,0.0)
\psdots[dotsize=0.12](10.780937,-1.0)
\usefont{T1}{ptm}{m}{n}
\rput(11.562344,0.81){$\mathcal{H}^{\{b_{3}\}}$}
\end{pspicture}
}

\end{exa}

\begin{rem}
We have the following toric interpretation of the Ehrhart polynomial $f_{\mathcal{B}}(m)$ (see, for example, \cite{HNPCayley}).
If $P_{i}$ is the convex hull of the origin and $b_{i}$ in $N_{\mathbb{R}}$, then one can consider the Minkowski sum $Q = P_{1} + \cdots + P_{n}$ in $N_{\mathbb{R}}$. 
Each $P_{i}$ determines a line bundle $L_{i}$ on the $d$-dimensional toric variety $Y$ corresponding to the polytope $Q$ and the toric variety corresponding to the polytope $P_{\mathcal{B}}$ is the $(d + n - 1)$-dimensional projective bundle $\mathbb{P}_{Y}( L_{1} \oplus \cdots \oplus L_{n} )$, with corresponding line bundle $\mathcal{O}(1)$. The Ehrhart polynomial is given by $f_{\mathcal{B}}(m) = \dim H^{0}(\mathbb{P}( L_{1} \oplus \cdots \oplus L_{n} ), \mathcal{O}(1)^{\otimes m})$, for positive integers $m$ \cite{FulIntroduction}.
\end{rem}

We conclude the introduction with a brief outline of the contents of the paper and note that all varieties and orbifolds will be over $\mathbb{C}$.
In Section \ref{construction}, we present a local construction of hypertoric varieties and use it to explain Theorem \ref{bussey}. In Section \ref{combinatorics}, we give a combinatorial proof of Theorem \ref{newbe}.

\section{Orbifold Cohomology of Hypertoric Varieties}\label{construction}

The goal of this section is to explain Theorem \ref{bussey}, which describes the dimensions of the orbifold cohomology groups of a hypertoric variety.
We will first review the construction of Lawrence toric varieties and refer the reader to Section 4 in \cite{HSToric}
for an equivalent presentation using Gale duality. 
We continue with the notation of the introduction and let
$H_{i, +} = \{ u \in M_{\mathbb{R}} \mid \langle u, b_{i} \rangle \geq r_{i}     \}$ and $H_{i, -} = \{ u \in M_{\mathbb{R}} \mid \langle u, b_{i} \rangle \leq r_{i}     \}$,
for $i = 1, \ldots, n$.
The hyperplane arrangement $\mathcal{H}$ divides $M_{\mathbb{R}}$ into locally closed \emph{cells} as follows: an element $u$ in $M_{\mathbb{R}}$ lies in the cell
given by the intersection of $\cap_{u \in H_{i,+}} H_{i,+} $ and $\cap_{u \in H_{i,-}} H_{i,-}$.
If $\sigma$ denotes the cone over the Lawrence polytope $P_{\mathcal{B}}$ in $N_{\mathbb{R}} \times \mathbb{R}^{n}$, then
$\mathcal{H}$ determines a simplicial fan refinement $\Sigma_{\mathcal{H}}$ of $\sigma$ \cite{HSToric}.
 If $C$ is a cell of $\mathcal{H}$, then there is a corresponding simplicial cone $\sigma_{C}$ in $\Sigma_{\mathcal{H}}$ with codimension equal to $\dim C$ and rays passing through the primitive integer vectors
$\{ (b_{i}, e_{i}) \mid C \subseteq H_{i,-} \}$ and $\{ (0, e_{i}) \mid C \subseteq H_{i,+} \}$.
 The cones in $\Sigma_{\mathcal{H}}$ not contained in the boundary of $\sigma$ are in bijection with the bounded cells of $\mathcal{H}$
 \cite[Theorem 4.7]{HSToric}.
In particular, the vertices of $\mathcal{H}$ correspond to the maximal cones of $\Sigma_{\mathcal{H}}$. The simplicial toric variety $X_{\mathcal{H}} = X(\Sigma_{\mathcal{H}})$ is the \emph{Lawrence toric variety} associated to $\mathcal{H}$. 
The proper torus-invariant subvarietes $V_{C}$ of $X_{\mathcal{H}}$ are in bijective correspondence with the bounded cells $C$ of $\mathcal{H}$.
The union of the proper torus-invariant subvarieties of $X_{\mathcal{H}}$ is called the \emph{core} of $X_{\mathcal{H}}$.

Fix a vertex $\gamma$ of $\mathcal{H}$ corresponding to a maximal element $F$ in $M_{\mathcal{B}}$ and to the 
maximal cone $\sigma_{\gamma}$ in $\Sigma_{\mathcal{H}}$.
We will describe the corresponding open toric subvariety $U_{\gamma}$ of $X_{\mathcal{H}}$. 
If $N_{F}$ denotes the sublattice of $N$ generated by $\{ b_{i} \mid b_{i} \in F \}$, then the elements of the finite group $N(F) = N/N_{F}$ are in bijective correspondence with  the elements of $\prod_{G \subseteq F} \BOX(G)\cap N$. 
Consider $\mathbb{A}^{d + n}$ with co-ordinates $\{ z_{i} \mid \gamma \in H_{i, -} \}$ and $\{ w_{i} \mid \gamma \in H_{i, +} \}$ and consider the embedding
$\iota: N(F) \hookrightarrow (\mathbb{C}^{*})^{d + n}$ such that if $g$ in $N(F)$ corresponds to
$v = \sum_{b_{i} \in G} \alpha_{i} b_{i} \in \BOX(G) \cap N$ with  $0 < \alpha_{i} < 1$, then,
for every
$b_{i}$ in $G$, $\iota(g)$  equals $e^{2 \pi i \alpha_{i}}$ in the co-ordinate corresponding to $z_{i}$ and equals $e^{2 \pi i (1 - \alpha_{i})}$ in the co-ordinate corresponding to $w_{i}$, and $\iota(g)$ equals $1$ in all other co-ordinates. 
The action of $(\mathbb{C}^{*})^{d + n}$ on $\mathbb{A}^{d + n}$ induces an
action of $N(F)$ on $\mathbb{A}^{d + n}$ such that the age of $g$ (see Subsection 7.1 \cite{AGVAlgebraic}) is given by
$\age(g) = \sum_{b_{i} \in G} \alpha_{i} + (1 - \alpha_{i}) = \dim G$. The corresponding open toric subvariety of $X_{\mathcal{H}}$ is given by  $U_{\gamma} = \mathbb{A}^{d + n}/ N(F)$ and the orbifold structure of $X_{\mathcal{H}}$ is locally induced by this quotient.
If $C$ is a bounded cell of $\mathcal{H}$, then the corresponding proper torus-invariant subvariety $V_{C}$ of $X_{\mathcal{H}}$ has $\dim V_{C} = \dim C$.
If $\gamma$ is not in the closure of $C$, then $U_{\gamma} \cap V_{C} = \emptyset$. If $\gamma$ lies in the closure of $C$, then $U_{\gamma} \cap V_{C}$ is defined by setting the co-ordinates $\{ z_{i} \mid C \subseteq H_{i,-} \}$ and $\{ w_{i} \mid C \subseteq H_{i, +} \}$ equal to zero.


The hypertoric variety $\mathcal{M}_{\mathcal{H}}$ associated to $\mathcal{H}$ is a $2d$-dimensional complete intersection in $X_{\mathcal{H}}$ with a canonical orbifold structure. We will describe $\mathcal{M}_{\mathcal{H}} \cap U_{\gamma}$ (c.f. \cite[Proposition 4.3]{JTOrbifold}).
If $b_{i} \notin F$, then we may write $b_{i} = \sum_{b_{j} \in F} a_{i,j} b_{j}$, for unique rational numbers $a_{i,j}$.
The ideal $I$ generated by $\langle z_{i} = \sum_{b_{j} \in F} a_{i,j}z_{j}w_{j} \mid 
u \in H_{i, -} \smallsetminus H_{i, +}  \rangle$ and $\langle w_{i} = \sum_{b_{j} \in F} a_{i,j}z_{j}w_{j} \mid 
u \in H_{i, +} \smallsetminus H_{i, -}  \rangle$,  defines a subvariety $V(I) \subseteq \mathbb{A}^{d + n}$ which is invariant under the action of
$N(F)$. The induced subvariety $V(I)/N(F)$ of $U_{\gamma}$ is equal to $\mathcal{M}_{\mathcal{H}} \cap U_{\gamma}$ and the orbifold structure of $\mathcal{M}_{\mathcal{H}}$ is locally induced by this quotient. It can be seen from this description that the core of $X_{\mathcal{H}}$ is contained in $\mathcal{M}_{\mathcal{H}}$.
If $\mathbb{A}^{2d}$ has co-ordinates $\{ z_{i} \mid b_{i} \in F \}$ and $\{ w_{i} \mid b_{i} \in F \}$, then the corresponding projection $\pi: \mathbb{A}^{d + n} \rightarrow \mathbb{A}^{2d}$ induces an isomorphism
$\pi: V(I) \rightarrow \mathbb{A}^{2d}$. The embedding $\iota: N(F) \hookrightarrow (\mathbb{C}^{*})^{d + n}$ factors as an embedding $\nu: N(F) \hookrightarrow (\mathbb{C}^{*})^{2d}$ followed by the inclusion of $(\mathbb{C}^{*})^{2d}$ into $(\mathbb{C}^{*})^{d + n}$ by adding $1$'s in the extra co-ordinates. The embedding $\nu$ induces an action of $N(F)$ on $\mathbb{A}^{2d}$ such that the projection $\pi$ is $N(F)$-equivariant and induces an isomorphism of orbifolds $\mathcal{M}_{\mathcal{H}} \cap U_{\gamma} \cong \mathbb{A}^{2d}/N(F)$.


We have seen in the introduction that any $G$ in $M_{\mathcal{B}}$ determines a hyperplane arrangement $\mathcal{H}^{G}$ with co-normal vectors $\mathcal{B}_{G}$. The corresponding Lawrence toric variety $X_{G}$ and hypertoric variety $\mathcal{M}_{G}$ may be regarded as subvarieties of $X_{\mathcal{H}}$ and $\mathcal{M}_{\mathcal{H}}$ respectively. For $g$ in $N(F)$
corresponding to
$v \in \BOX(G) \cap N$,
one verifies that,  as varieties,  $X_{G} \cap U_{\gamma} = (\mathbb{A}^{d + n})^{g}/N(F)$
and $\mathcal{M}_{G} \cap U_{\gamma} \cong (\mathbb{A}^{2d})^{g}/N(F)$, where $(\mathbb{A}^{d + n})^{g}$ and  $(\mathbb{A}^{2d})^{g}$ are the subvarieties of $(\mathbb{A}^{d + n})$ and $(\mathbb{A}^{2d})$ respectively that are invariant under the action of $g$.

\begin{exa}
Consider the notation of Example \ref{visit} and fix the vertex $\gamma = (-1/2, 0)$ of $\mathcal{H}$ corresponding to $F = \{ b_{2}b_{4} \}$ in $M_{\mathcal{B}}$.
The open toric subvariety $U_{\gamma}$ of $X_{\mathcal{H}}$ is the quotient of $\mathbb{A}^{6}$ with co-ordinates $\{z_{1},z_{2},z_{3},
z_{4},w_{2},w_{4}\}$ by the action of $N(F) = \mathbb{Z}/2\mathbb{Z}$, acting via multiplication by $(1,-1,1,-1,-1,-1)$.
The open subvariety $\mathcal{M}_{\mathcal{H}} \cap U_{\gamma}$ of the hypertoric variety $\mathcal{M}_{\mathcal{H}}$ is the complete intersection defined by $z_{1} = \frac{1}{2} z_{2}w_{2} + \frac{1}{2} z_{4}w_{4}$ and $z_{3} = -z_{2}w_{2} - z_{4}w_{4}$.
The core of $X_{\mathcal{H}}$ consists of two weighted projective spaces $\mathbb{P}(1,1,2)$ that intersect along their unique singular points at $V_{\gamma} = X_{F} = \mathcal{M}_{F}$.
\end{exa}

\begin{rem}
Hausel and Sturmfels proved that the embedding of the core of $X_{\mathcal{H}}$ into either $X_{\mathcal{H}}$ or $\mathcal{M}_{\mathcal{H}}$ induces an isomorphism of cohomology rings over $\mathbb{Z}$ \cite[Lemma 6.5]{HSToric} and used this to prove Theorem \ref{hussey} and Theorem \ref{rubbish}.
\end{rem}


Orbifold cohomology was introduced by Chen and Ruan in \cite{CRNew} and associates to an orbifold $Y$ a graded $\mathbb{Q}$-algebra $H^{*}_{\orb}(Y, \mathbb{Q})$. More specifically, one associates to $Y$ an orbifold $I(Y)$, called the \emph{inertia orbifold} of $Y$, such that, 
if $Y = X/G$ for some smooth variety $X$ and finite abelian group $G$, then $I(Y) = \coprod_{g \in G} X^{g}/G$, where $X^{g}$ is the subvariety of $X$ fixed by $g$. If $Y_{1}, \ldots, Y_{r}$ denote the connected components of $I(Y)$, then,
for any $i \in \mathbb{Q}$, Chen and Ruan defined the $i$'th 
orbifold cohomology group of $Y$ by
\begin{equation*}
H_{\orb}^{i}(Y, \mathbb{Q}) = \bigoplus_{j = 1}^{r} H^{i - 2\age(Y_{j})}(Y_{j}, \mathbb{Q}),
\end{equation*}
where $\age(Y_{j})$ is the \emph{age} of $Y_{j}$.
We have $I(U_{\gamma}) = \coprod_{G \subseteq F} \coprod_{  v \in \BOX(G) \cap N } (\mathbb{A}^{d + n})^{g}/N(F)$, $I(\mathcal{M}_{\mathcal{H}} \cap U_{\gamma}) = \coprod_{G \subseteq F} \coprod_{  v \in \BOX(G) \cap N } (\mathbb{A}^{2d})^{g}/N(F)$, and, in both cases, the age of the connected component corresponding to a pair $(v,G)$ is equal to $\dim G$.
In fact,  we have isomorphisms of varieties, $I(X_{\mathcal{H}}) \cong \coprod_{G \in M} \coprod_{  v \in \BOX(G) \cap N } X_{G}$ and
$I(\mathcal{M}_{\mathcal{H}}) \cong \coprod_{G \in M} \coprod_{  v \in \BOX(G) \cap N } \mathcal{M}_{G}$ \cite{JTNote}. We conclude that Theorem \ref{hussey} and the above computation of inertia orbifolds implies Theorem \ref{bussey}, which computes the orbifold cohomology of $X_{\mathcal{H}}$ and $\mathcal{M}_{\mathcal{H}}$.

\section{Ehrhart Theory for Lawrence Polytopes}\label{combinatorics}

The goal of this section is to give a combinatorial proof of Theorem \ref{newbe}. We will first show that
\begin{equation}\label{done}
\delta_{\mathcal{B}}(t) =
\sum_{F \in M} \# \left( \BOX(F) \cap N \right) t^{\dim F} h^{bd}_{F}(t).
\end{equation}
The proof should be compared with the combinatorial proof of Theorem \ref{oldie} (see \cite{YoWeightI}) and the proof of \cite[Theorem 1.1]{PayEhrhart}.  Recall that $\sigma$ denotes the cone over $P_{\mathcal{B}}$ in $N_{\mathbb{R}} \times \mathbb{R}^{n}$ and that $\Sigma_{\mathcal{H}}$ is a fan refining $\sigma$, such that  the cones $\sigma_{C}$ in $\Sigma_{\mathcal{H}}$ that do not lie in the boundary of $\sigma$ are in bijection with the bounded cells $C$ of $\mathcal{H}$. The cone $\sigma_{C}$ has codimension equal to $\dim C$ and its rays correspond to the primitive integer vectors
$\{ (b_{i}, e_{i}) \mid C \subseteq H_{i,-} \}$ and $\{ (0, e_{i}) \mid C \subseteq H_{i,+} \}$. If $F$ is an element of $M_{\mathcal{B}}$ and $w = \sum_{b_{i} \in F} a_{i}b_{i} \in \BOX(F) \cap N$, for some $0< a_{i} < 1$, then let
$l(w) = \sum_{b_{i} \in F} (a_{i}b_{i},e_{i}) \in \sigma \cap (N \times \mathbb{Z}^{n})$. 
If $v$ is a lattice point in the interior of $\sigma$ and $C$ is the largest bounded cell in $\mathcal{H}$ such that $v \in \sigma_{C}$, then we have a unique expression
\[
v = \: l(w) + \sum_{b_{i} \notin F, C \subseteq H_{i,-}} (b_{i}, e_{i}) + \sum_{b_{i} \notin F, C \subseteq H_{i,+} } (0, e_{i})
+ \sum_{C \subseteq H_{i,-}} \alpha_{i} (b_{i}, e_{i}) + \sum_{C \subseteq H_{i,+}} \beta_{i}(0, e_{i}),
\]
where $w \in \BOX(F) \cap N$, for some $F$ in $M_{\mathcal{B}}$ such that $C$ is a bounded cell in $\mathcal{H}^{F}$, and $\alpha_{i}, \beta_{i}$ are non-negative integers. Conversely, given $w \in \BOX(F) \cap N$, a bounded cell $C$ in $\mathcal{H}^{F}$
and non-negative integers $\{ \alpha_{i} \mid C \subseteq H_{i,-} \}$ and $\{ \beta_{i} \mid C \subseteq H_{i,+} \}$,
the lattice point
$v = l(w) + \sum_{b_{i} \notin F, C \subseteq H_{i,-}} (b_{i}, e_{i}) + \sum_{b_{i} \notin F, C \subseteq H_{i,+} } (0, e_{i})
+ \sum_{C \subseteq H_{i,-}} \alpha_{i} (b_{i}, e_{i}) + \sum_{C \subseteq H_{i,+}} \beta_{i}(0, e_{i})$, lies in the interior of $\sigma$ and $C$ is the largest bounded cell in $\mathcal{H}$ such that $v \in \sigma_{C}$. Recall from the proof of Theorem \ref{punter} that
$\psi: N \times \mathbb{Z}^{n} \rightarrow \mathbb{Z}$ is given by $\psi(\lambda + \sum_{i = 1}^{n} \mu_{i}e_{i} ) = \sum_{i = 1}^{n} \mu_{i}$, for $\lambda \in N$ and integers $\mu_{i}$, and that $\psi^{-1}(m) \cap \sigma = mP_{\mathcal{B}}$ and $\psi(l(w)) = \dim F$. It is a standard result of Ehrhart theory (see, for example, \cite{BRComputing}) that
\[ \sum_{m \geq 1} \# \left( \Int(mP_{\mathcal{B}}) \cap (N \times \mathbb{Z}^{n}) \right)t^{m} = \frac{t^{d + n} \delta_{\mathcal{B}}(t^{-1})}{(1 - t)^{d + n}}. \]
Using the above facts and setting $\mathcal{H}_{\bd}^{F}$ to be the collection of bounded cells in $\mathcal{H}^{F}$, we calculate
\begin{align*}
t^{d + n} \delta_{\mathcal{B}}(t^{-1}) &= (1 - t)^{d + n} \sum_{ v \in \Int(\sigma) \cap (N \times \mathbb{Z}^{n})} t^{\psi(v)} \\
&= (1 - t)^{d + n} \sum_{F \in M} \# \left( \BOX(F) \cap N \right) t^{\dim F} \sum_{C \in \mathcal{H}_{\bd}^{F}} \frac{t^{d + n - \dim C -2\dim F}}{( 1- t)^{d + n - \dim C}} \\
&= t^{d + n} \sum_{F \in M} \# \left( \BOX(F) \cap N \right) t^{-\dim F} \sum_{C \in \mathcal{H}_{\bd}^{F}} (t^{-1} - 1)^{\dim C},
\end{align*}
and (\ref{done}) follows immediately as desired. We present below the remainder of the proof of Theorem \ref{newbe}.

\begin{proof}
Observe that the constant coefficient in $h_{F}(t)$ is $1$ and the coefficient of $t$ in $h_{\{ 0 \}}(t)$ is
$n - d$. The elements of $M_{\mathcal{B}}$ of dimension $1$ correspond to the elements $b_{i}$ in $\mathcal{B}$ and $\BOX(\{b_{i}\} ) \cap N$ consists of $a_{i} - 1$ lattice points. We deduce that $\delta_{1} = n - d + \sum_{i = 1}^{n} (a_{i} - 1) = \sum_{i = 1}^{n}a_{i} - d$. The leading term of
$h^{bd}_{F}(t)$ is $R^{bd}_{F} t^{\codim F}$ and hence $\delta_{d} = \sum_{F \in M} \# \left( \BOX(F) \cap N \right) R^{bd}_{F}$ and
$\delta_{d + 1} = \ldots = \delta_{d + n - 1} = 0$. It is a standard fact of Ehrhart theory (see, for example, \cite{BRComputing}) that the normalised volume of $P_{B}$ is equal to $\frac{1}{(d + n - 1)!} \delta_{\mathcal{B}}(1)$. The last statement follows from the observation that
$h_{F}(1) = h^{\bd}_{F}(1) = V_{F}$.
\end{proof}

\begin{rem}\label{smashing}
If $P_{1}, \ldots, P_{n}$ are lattice polytopes in $N_{\mathbb{R}}$ and $\{e_{1}, \ldots, e_{n} \}$ is the standard basis of $\mathbb{R}^{n}$, then the \emph{Cayley sum} $P = P_{1} \ast \cdots \ast P_{n}$ is the convex hull of
$P_{1} \times \{e_{1}\}$, \ldots, $P_{n} \times \{e_{n}\}$ in $N \times \mathbb{R}^{n}$. If the affine span of the union of the $P_{i}$ is $N_{\mathbb{R}}$, then $P$ is a $(d + n - 1)$-dimensional lattice polytope. Setting $P_{i}$ to be the convex hull of the origin and $b_{i}$, we
recover the Lawrence polytope $P_{\mathcal{B}}$.
The degree $s$ of a lattice polytope $Q$ is the degree of its Ehrhart $\delta$-polynomial and it is a standard fact that  $(\dim Q + 1 - s)Q$ is the smallest dilate of $Q$
that contains a lattice point in its relative interior (see, for example, \cite{BRComputing}). Cayley sums provide examples of lattice polytopes with small degree relative to their dimension since the degree of $P$ is at most $d$ \cite[Proposition]{BNMultiples}. In \cite{BNMultiples}, Batyrev and Nill 
used Lawrence polytopes in the classification of lattice polytopes of degree $1$.
\end{rem}

\begin{rem}
It is standard result of Ehrhart theory (see, for example, \cite{BRComputing}) that $\delta_{1} = f_{\mathcal{B}}(1) - \dim P_{\mathcal{B}} - 1$
and hence Theorem \ref{newbe} implies that $P_{\mathcal{B}}$ contains $\sum_{i = 1}^{n} (a_{i} + 1)$ lattice points. More specifically, the lattice points in $P_{\mathcal{B}}$ are
$\{ (\lambda_{i}v_{i}, e_{i}) \mid 0 \leq \lambda_{i} \leq a_{i}, 1 \leq i \leq n \}$. We noted in Remark \ref{smashing} that $mP_{\mathcal{B}}$ contains no interior lattice points for $1 \leq m \leq n - 1$. It is a standard fact that $\delta_{d} = \# ( \Int(nP_{\mathcal{B}}) \cap (N \times \mathbb{Z}^{n}) )$ (see, for example, \cite{BRComputing}), 
and hence Theorem \ref{newbe} implies that $nP$ contains $\sum_{F \in M} \# \left( \BOX(F) \cap N \right) R^{bd}_{F}$ interior lattice points. More specifically, if $w$ is a lattice point in $\BOX(F)$ and $R$ is a bounded region in $\mathcal{H}^{F}$, then the corresponding lattice point in the relative interior of $nP_{\mathcal{B}}$ is $l(w) + \sum_{b_{i} \notin F, R \subseteq H_{i,-}} (b_{i}, e_{i}) + \sum_{b_{i} \notin F,  R \subseteq H_{i,+}} (0, e_{i}) \in \sigma_{R}$.
\end{rem}

\begin{rem}
We present an alternative combinatorial proof of the formula
\begin{equation*}
\delta_{\mathcal{B}}(t) =
\sum_{F \in M} \# \left( \BOX(F) \cap N \right) t^{\dim F} h_{F}(t).
\end{equation*}
If $v$ is a lattice point in $\sigma$, then $v$ lies in some maximal cone $\sigma_{\gamma}$ of $\Sigma_{\mathcal{H}}$, and we may write
$v = \sum_{\gamma \in H_{i,-}} \alpha_{i} (b_{i}, e_{i}) + \sum_{\gamma \in H_{i,+}} \beta_{i} (0, e_{i})$, for unique non-negative integers $\alpha_{i}$ and $\beta_{i}$. The element $G$ of $M_{\mathcal{B}}$ generated by $\{ b_{i} \mid \alpha_{i}, \beta_{i} \neq 0 \}$ and the sets
 $I_{-} = \{ i \mid \alpha_{i} \neq 0, b_{i} \notin \Span G  \}$ and
$I_{+} = \{ i \mid \beta_{i} \neq 0, b_{i} \notin \Span G  \}$
are independent of the choice of maximal cone $\sigma_{\gamma}$.
Note that $b_{i} \in \Span G$ if and only if  $M_{G,\mathbb{R}} = \cap_{b_{i} \in G} H_{i}$ is contained in either  $H_{i,-}$ or $H_{i,+}$.
If we consider the projection $\phi_{G}: N \rightarrow  N/ (N \cap \Span G)$ and the hyperplane arrangement $\{ \phi_{G}(H_{i}) \mid i \in I = I_{-} \cup I_{+}  \}$, then the intersection of the half spaces $\{ \phi_{G}(H_{i, -}) \mid  i \in I_{-} \}$ and $\{ \phi_{G}(H_{i, +}) \mid  i \in I_{+} \}$ is a region in the hyperplane arrangement and, moreover, every region in the hyperplane arrangement has this form.
We deduce that $v$ has a unique expression
\[
v = \: l(w) + \sum_{b_{i} \in G \smallsetminus F} ((b_{i}, e_{i}) + (0, e_{i}))
+ \sum_{M_{G,\mathbb{R}} \subseteq H_{i,-}} \mu_{i} (b_{i}, e_{i}) + \sum_{M_{G,\mathbb{R}} \subseteq H_{i,+}} \mu_{i}'(0, e_{i})
+ \sum_{i \in I_{-}} \nu_{i}(b_{i}, e_{i}) + \sum_{i \in I_{+}} \nu_{i}'(0, e_{i}),
\]
where $w \in \BOX(F) \cap N$, for some $F \subseteq G$, $\mu_{i}, \mu_{i}'$ are non-negative integers and $\nu$, $\nu'$ are positive integers.
Conversely, consider elements $F \subseteq G$ in $M_{\mathcal{B}}$, a subset $I$ in $\{ i \mid b_{i} \notin \Span G \}$ and a region $R$ in the hyperplane arrangement $\{ \phi_{G}(H_{i}) \mid i \in I  \}$. If we set $I_{-} = \{ i \mid R \subseteq \phi_{G}(H_{i, -}) \}$ and $I_{+} = \{ i \mid R \subseteq \phi_{G}(H_{i, +}) \}$, and consider some $w \in \BOX(F) \cap N$, some non-negative integers $\mu_{i}, \mu_{i}'$ and some positive integers $\nu$, $\nu'$, then the right hand side of the above expression gives a lattice point $v$ in $\sigma$.
By a theorem of Zaslavsky \cite{ZasFacing}, the number of regions $r_{G,I}$ in the hyperplane arrangement $\{ \phi_{G}(H_{i}) \mid i \in I  \}$
is equal to the number of elements $G'$ in $M_{\mathcal{B}}$ satisfying $G \subseteq G'\subseteq G \cup \{ b_{i} \mid i \in I \}$.
Setting
$I_{G} = \{ i \mid b_{i} \notin \Span G \}$, we compute
\begin{align*}
\delta_{\mathcal{B}}(t) &= (1 - t)^{d + n} \sum_{ v \in \sigma \cap (N \times \mathbb{Z}^{n})} t^{\psi(v)} \\
&= (1 - t)^{d + n} \sum_{F \in M} \# \left( \BOX(F) \cap N \right) t^{\dim F} \sum_{F \subseteq G}  \frac{t^{2(\dim G - \dim F)}}{(1 - t)^{\dim G + n - |I_{G}|}} \sum_{I \subseteq I_{G} } r_{G,I} \frac{t^{|I|}}{(1 - t)^{|I| }} \\
&= \sum_{F \in M} \# \left( \BOX(F) \cap N \right) t^{-\dim F} \sum_{F \subseteq G} t^{2\dim G} (1 - t)^{d - \dim G + |I_{G}|} \sum_{G \subseteq G'} \frac{t^{\dim G' - \dim G}}{(1 - t)^{|I_{G}|}} \\
&= \sum_{F \in M} \# \left( \BOX(F) \cap N \right) t^{\dim F }  \sum_{F \subseteq G'} t^{\dim G' - \dim F }(1 - t)^{d - \dim G'} \sum_{F \subseteq G \subseteq G'} t^{\dim G - \dim F} (1 - t)^{\dim G' - \dim G} \\
&= \sum_{F \in M} \# \left( \BOX(F) \cap N \right) t^{\dim F} h_{F}(t).
\end{align*}


\end{rem}

\bibliographystyle{amsplain}
\bibliography{alan}

\def\cprime{$'$}
\providecommand{\bysame}{\leavevmode\hbox to3em{\hrulefill}\thinspace}
\providecommand{\MR}{\relax\ifhmode\unskip\space\fi MR }
\providecommand{\MRhref}[2]{%
  \href{http://www.ams.org/mathscinet-getitem?mr=#1}{#2}
}
\providecommand{\href}[2]{#2}
\begin{thebibliography}{10}

\bibitem{AGVAlgebraic}
Dan Abramovich, Tom Graber, and Angelo Vistoli, \emph{Algebraic orbifold
  quantum products}, Orbifolds in mathematics and physics (Madison, WI, 2001),
  Contemp. Math., vol. 310, Amer. Math. Soc., Providence, RI, 2002, pp.~1--24.
  \MR{MR1950940 (2004c:14104)}

\bibitem{BNMultiples}
Victor Batyrev and Benjamin Nill, \emph{Multiples of lattice polytopes without
  interior lattice points}, Mosc. Math. J. \textbf{7} (2007), no.~2, 195--207,
  349. \MR{MR2337878 (2008g:52018)}

\bibitem{BSLawrence}
Margaret Bayer and Bernd Sturmfels, \emph{Lawrence polytopes}, Canad. J. Math.
  \textbf{42} (1990), no.~1, 62--79. \MR{MR1043511 (91e:52023)}

\bibitem{BRComputing}
Matthias Beck and Sinai Robins, \emph{Computing the continuous discretely},
  Undergraduate Texts in Mathematics, Springer, New York, 2007, Integer-point
  enumeration in polyhedra. \MR{MR2271992 (2007h:11119)}

\bibitem{BDGeometry}
Roger Bielawski and Andrew~S. Dancer, \emph{The geometry and topology of toric
  hyperk\"ahler manifolds}, Comm. Anal. Geom. \textbf{8} (2000), no.~4,
  727--760. \MR{MR1792372 (2002c:53078)}

\bibitem{CROrbifold}
Weimin Chen and Yongbin Ruan, \emph{Orbifold {G}romov-{W}itten theory},
  Orbifolds in mathematics and physics (Madison, WI, 2001), Contemp. Math.,
  vol. 310, Amer. Math. Soc., Providence, RI, 2002, pp.~25--85. \MR{MR1950941
  (2004k:53145)}

\bibitem{CRNew}
\bysame, \emph{A new cohomology theory of orbifold}, Comm. Math. Phys.
  \textbf{248} (2004), no.~1, 1--31. \MR{MR2104605 (2005j:57036)}

\bibitem{HNPCayley}
Hasse Christian, Nill Benjamin, and Payne Sam, \emph{Cayley decompositions of
  lattice polytopes and upper bounds for $h^{*}$-polynomials}, arXiv:0804.3667,
  2008.

\bibitem{ehrhartpolynomial}
Eug{\`e}ne Ehrhart, \emph{Sur les poly\`edres rationnels homoth\'etiques \`a
  {$n$}\ dimensions}, C. R. Acad. Sci. Paris \textbf{254} (1962), 616--618.
  \MR{24 \#A714}

\bibitem{FulIntroduction}
William Fulton, \emph{Introduction to toric varieties}, Annals of Mathematics
  Studies, vol. 131, Princeton University Press, Princeton, NJ, 1993, , The
  William H. Roever Lectures in Geometry. \MR{MR1234037 (94g:14028)}

\bibitem{GHOrbifold}
Rebecca~F. Goldin and Megumi Harada, \emph{Orbifold cohomology of hypertoric
  varieties}, arXiv:math/0607421, 2006.

\bibitem{GruConvex}
Branko Gr{\"u}nbaum, \emph{Convex polytopes}, second ed., Graduate Texts in
  Mathematics, vol. 221, Springer-Verlag, New York, 2003, Prepared and with a
  preface by Volker Kaibel, Victor Klee and G\"unter M.\ Ziegler. \MR{MR1976856
  (2004b:52001)}

\bibitem{HPProperties}
Megumi Harada and Nicholas Proudfoot, \emph{Properties of the residual circle
  action on a hypertoric variety}, Pacific J. Math. \textbf{214} (2004), no.~2,
  263--284. \MR{MR2042933 (2004k:53139)}

\bibitem{HauQuaternionic}
Tam{\'a}s Hausel, \emph{Quaternionic geometry of matroids}, Cent. Eur. J. Math.
  \textbf{3} (2005), no.~1, 26--38 (electronic). \MR{MR2110782 (2005m:53070)}

\bibitem{HSToric}
Tam{\'a}s Hausel and Bernd Sturmfels, \emph{Toric hyper{K}\"ahler varieties},
  Doc. Math. \textbf{7} (2002), 495--534 (electronic). \MR{MR2015052
  (2004i:53054)}

\bibitem{HibSome}
Takayuki Hibi, \emph{Some results on {E}hrhart polynomials of convex
  polytopes}, Discrete Math. \textbf{83} (1990), no.~1, 119--121. \MR{MR1065691
  (91g:52008)}

\bibitem{HibEhrhart}
\bysame, \emph{Ehrhart polynomials of convex polytopes, {$h$}-vectors of
  simplicial complexes, and nonsingular projective toric varieties}, Discrete
  and computational geometry (New Brunswick, NJ, 1989/1990), DIMACS Ser.
  Discrete Math. Theoret. Comput. Sci., vol.~6, Amer. Math. Soc., Providence,
  RI, 1991, pp.~165--177. \MR{MR1143294 (92j:52018)}

\bibitem{HibLower}
\bysame, \emph{A lower bound theorem for {E}hrhart polynomials of convex
  polytopes}, Adv. Math. \textbf{105} (1994), no.~2, 162--165. \MR{MR1275662
  (95b:52018)}

\bibitem{HibStar}
\bysame, \emph{Star-shaped complexes and {E}hrhart polynomials}, Proc. Amer.
  Math. Soc. \textbf{123} (1995), no.~3, 723--726. \MR{MR1249883 (95d:52012)}

\bibitem{JTNote}
Yunfeng Jiang and Hsian-Hua Tseng, \emph{Note on orbifold chow ring of
  semi-projective toric deligne-mumford stacks}, arXiv:math/0606322, 2006.

\bibitem{JTOrbifold}
\bysame, \emph{The orbifold chow ring of hypertoric deligne-mumford stacks},
  arXiv:math/0512199, 2005.

\bibitem{KonCohomology}
Hiroshi Konno, \emph{Cohomology rings of toric hyperk\"ahler manifolds},
  Internat. J. Math. \textbf{11} (2000), no.~8, 1001--1026. \MR{MR1797675
  (2001k:53089)}

\bibitem{MnevUniversality}
N.~E. Mn{\"e}v, \emph{The universality theorems on the classification problem
  of configuration varieties and convex polytopes varieties}, Topology and
  geometry---Rohlin Seminar, Lecture Notes in Math., vol. 1346, Springer,
  Berlin, 1988, pp.~527--543. \MR{MR970093 (90a:52013)}

\bibitem{ProSurvey}
Proudfoot Nicholas, \emph{A survey of hypertoric geometry and topology},
  arXiv:0705.4236, 2007.

\bibitem{PayEhrhart}
Sam Payne, \emph{Ehrhart series and lattice triangulations},
  arXiv:math/0702052, to appear in Discr. Comput. Geom, 2007.

\bibitem{SanTriangulations}
Francisco Santos, \emph{Triangulations of oriented matroids}, Mem. Amer. Math.
  Soc. \textbf{156} (2002), no.~741, viii+80. \MR{MR1880595 (2003b:52013)}

\bibitem{StaHilbert2}
Richard~P. Stanley, \emph{Hilbert functions of graded algebras}, Advances in
  Math. \textbf{28} (1978), no.~1, 57--83. \MR{MR0485835 (58 \#5637)}

\bibitem{StaDecompositions}
\bysame, \emph{Decompositions of rational convex polytopes}, Ann. Discrete
  Math. \textbf{6} (1980), 333--342, Combinatorial mathematics, optimal designs
  and their applications (Proc. Sympos. Combin. Math. and Optimal Design,
  Colorado State Univ., Fort Collins, Colo., 1978). \MR{MR593545 (82a:52007)}

\bibitem{StaHilbert1}
\bysame, \emph{On the {H}ilbert function of a graded {C}ohen-{M}acaulay
  domain}, J. Pure Appl. Algebra \textbf{73} (1991), no.~3, 307--314.
  \MR{MR1124790 (92f:13017)}

\bibitem{StaMonotonicity}
\bysame, \emph{A monotonicity property of {$h$}-vectors and {$h\sp
  *$}-vectors}, European J. Combin. \textbf{14} (1993), no.~3, 251--258.
  \MR{MR1215335 (94f:52016)}

\bibitem{YoWeightI}
Alan Stapledon, \emph{Weighted ehrhart theory and orbifold cohomology}, Adv.
  Math. (2008), doi:10.1016/j.aim.2008.04.010.

\bibitem{SwaElements}
Ed~Swartz, \emph{{$g$}-elements of matroid complexes}, J. Combin. Theory Ser. B
  \textbf{88} (2003), no.~2, 369--375. \MR{MR1983365 (2004g:05041)}

\bibitem{ZasFacing}
Thomas Zaslavsky, \emph{Facing up to arrangements: face-count formulas for
  partitions of space by hyperplanes}, Mem. Amer. Math. Soc. \textbf{1} (1975),
  no.~issue 1, 154, vii+102. \MR{MR0357135 (50 \#9603)}

\bibitem{ZieNon}
G{\"u}nter~M. Ziegler, \emph{Non-rational configurations, polytopes, and
  surfaces}, arXiv:0710.4453, 2007.

\bibitem{ZieLectures}
\bysame, \emph{Lectures on polytopes}, Graduate Texts in Mathematics, vol. 152,
  Springer-Verlag, New York, 1995. \MR{MR1311028 (96a:52011)}

\end{thebibliography}

\end{document}